\documentclass[preprint,12pt]{elsarticle}

\usepackage{amssymb}
\usepackage{stmaryrd}
\usepackage{tipa}
\usepackage{amsfonts,amsmath,latexsym}
\usepackage{mathrsfs}  
\usepackage{amsbsy}
\usepackage{indentfirst}
\usepackage{amsmath}
\usepackage{amsfonts}
\usepackage{bbm}
\usepackage{mathrsfs}
\usepackage{dsfont}
\usepackage{graphicx}

\numberwithin{equation}{section}

\newtheorem{theorem}{\textbf{Theorem}}[section]

\newtheorem{definition}{\textbf{Definition}}

\newtheorem{remark}[theorem]{Remark}

\newcommand{\beqnar}{\begin{eqnarray*}}
\newcommand{\eeqnar}{\end{eqnarray*}}
\newcommand{\ba}{\begin{array}}
\newcommand{\ea}{\end{array}}

\journal{}

\allowdisplaybreaks[4]

\begin{document}

\begin{frontmatter}



\title{On Catoni's M-Estimation}


 \author{Pengtao Li}

\address{ School of Mathematics, Shandong University, Jinan, 250100,  China
    }

\author{Hanchao Wang\footnote{Corresponding author, email: wanghanchao@sdu.edu.cn.}}

\address{ Institute for Financial Studies,  Shandong University,  Jinan,  250100, China
    }

 \begin{abstract}
Catoni proposed a robust M-estimator and gave the deviation inequality for one fixed test function. The present paper is devoted to the uniform concentration inequality for a family of test functions. As an application, we consider empirical risk minimization for heavy-tailed losses.

\end{abstract}


\end{frontmatter}



\section{Introduction and main results}

Let ($\Omega,\mathcal{F},\mathbb{P} $) be a probability triple on which a sequence of independent and identically distributed (i.i.d.) random variables $(X_i)_{i=1}^{n}$ is defined and $X$ be independent from $(X_i)_{i=1}^{n}$ with the same marginal distribution $F$. Estimating the expected value $\mu$ of $X$ based on the sample $(X_i)_{i=1}^{n}$ is one of the most fundamental problems in statistics. 

The most natural choice is the empirical mean $$\mu_n=\frac{1}{n}\displaystyle\sum_{i=1}^n X_i.$$
For Gaussian distribution, the empirical mean $\mu_n$ has optimal minimax mean square error (see Catoni \cite{c} for details). However, when the distribution is heavy-tailed, meaning that only finite moments of low order exist (in the context of this paper, "low order" will usually mean the range between 2 and 3), the empirical mean $\mu_n$ is far from optimal (see Lugosi and Mendelson \cite{l1} or \cite{l2} for details). Catoni \cite{c} proposed an M-estimator, which is called Catoni's estimator, to cope with the heavy-tailed data.

Some notations are needed to be introduced at the beginning. Let $\Psi$ be a set of real-valued functions. For each $f\in\Psi$, denote the expectation $\mathbb{E}f(X)$ by $m_f$. For many problems in statistical learning theory, it is required to estimate $m_f$ for all $f\in\Psi$ simultaneously. For instance, in the maximum likelihood estimation context, $\Psi=\{\log p_{\theta}(\cdot),\theta\in\boldsymbol{\Theta}\}$ is a family of probability density function with respect to a $\sigma-$finite measure $\tau$ and $\frac{d\mathbb{P}}{d\tau}=p_{\theta^*}$ for $\theta^*\in\boldsymbol{\Theta}$.

Inspired by Catoni \cite{c}, for any given $f\in\Psi$, any non-decreasing differentiable real-valued function $\rho:\mathbb{R}\rightarrow\mathbb{R}$ and all $\mu\in\mathbb{R}$, the quantity $\hat{r}_{f,\rho}(\mu)$ is defined as follows,
$$\hat{r}_{f,\rho}(\mu)=\frac{1}{n\alpha}\displaystyle\sum_{i=1}^{n}\rho\left(\alpha(f(X_i)-\mu)\right),$$
where $\alpha>0$ is a parameter to be optimized. The value $\hat{\mu}_{f,\rho}$ for which $\hat{r}_{f,\rho}(\hat{\mu}_{f,\rho})=0$ is of our interest (If the value is not unique, we can pick one of them arbitrarily).

When we choose $\rho_0(x)=x$, $\hat{\mu}_{f,\rho_0}$ is virtually the empirical mean $\mu_n$.

To obtain Catoni's estimator, introduce the non-decreasing differentiable truncation function $\phi:\mathbb{R}\rightarrow\mathbb{R}$ such that
\begin{equation}\label{catoni}
    -\log(1-x+\frac{x^2}{2})\leq\phi(x)\leq\log(1+x+\frac{x^2}{2}).
\end{equation}

Let $\rho(x)=\phi(x)$, and then we get Catoni's estimator $\hat{\mu}_{f,\phi}$. Since $\phi$ is any fixed function satisfying (\ref{catoni}), we denote Catoni's estimator by $\hat{\mu}_{f}$ for simplicity.

For any fixed $f\in\Psi$, assuming that $Var[f(X)]\leq \sigma^2$ for some $\sigma^2>0$, for $\delta\in(0,1)$ such that $n>2\log(1/\delta)$, with proper choice of $\alpha>0$,
Catoni \cite{c} showed that
\begin{equation}\label{deviation}
    \mathbb{P}\left(|m_f-\hat{\mu}_f|\geq\sqrt{\frac{2\sigma^2\log(1/\delta)}{n(1-(2\log(1/\delta)/n))}}\right)\leq2\delta. 
\end{equation}

Rewrite (\ref{deviation}), and we know that 
\begin{equation}\label{deviations}
\mathbb{P}\big(|m_f-\hat{\mu}_f|\geq x\big)\leq2\,\exp\left(-\frac{nx^2}{2(\sigma^2+x^2)}\right)\,\quad(x>0)
\end{equation}
under the assumption that $\displaystyle\sup_{f\in\Psi} Var[f(X)]\leq \sigma^2$ for some $\sigma^2>0$.

This gives Catoni's estimator for a single $f\in\Psi$ and its deviation inequality. 

As for robust mean estimators for heavy-tailed distributions, Catoni's estimator is not the only possible estimator. Median-of-means (MOM) estimator proposed by Nemirovsky and Yudin \cite{ny} is also a powerful robust estimator. The median-of-means estimate is obtained by dividing the data into several small blocks of roughly equal size, calculating the empirical mean within each block, and then taking the median of these obtained values. Lugosi and Mendelson \cite{l2} demonstrate the \textit{sub-Gaussian} performance of the MOM estimator under some conditions (see van der Vaart and Wellner \cite{vw} for more about ``sub-Gaussian''). Moreover, Minsker \cite{m}  constructed a family of estimators, which we call Minsker's estimators, and they can be viewed as a combination of Catoni's estimator and the MOM technique. For more robust mean estimators under heavy-tailed distributions, see \cite{l2}. Here, we mainly focus on Catoni's estimators.

With Catoni's estimator and deviation inequality for a single $f\in\Psi$ in hand, an ensuing problem is the uniform concentration inequality for $f\in\Psi$. Actually, in statistical learning theory, it is natural and vital to consider the suprema of the deviations for all the simultaneous mean estimations after obtaining all the mean estimators simultaneously. The theory of empirical process has extensively studied this problem, assuming either $f\in\Psi$ are uniformly bounded or $f(X)$ have sub-Gaussian tails for all $f\in\Psi$. In such cases, empirical means perform well. The interested reader is referred to Bartlett and Mendelson \cite{bar}, Koltchinskii \cite{kol} and van der Vaart and Wellner \cite{vw}. When $f\in\Psi$ are not uniformly bounded, and random variables $f(X)$ have a heavy tail, we replace empirical means with robust proxies of the expectation and consider uniform tail-bounds for the deviations of the simultaneous robust mean estimations. For example, for Minsker's estimators, Minkser \cite{m} gave the uniform bound.

The main contribution of this paper is providing uniform exponential inequalities for the deviations of the stochastic process defined by Catoni's estimators.

As in empirical process theory, some conditions on $\mathscr{X}=\{m_f-\hat{\mu}_f,\,f\in\Psi\}$ are needed to make this possible. The norm we use is the essential supremum norm
\begin{equation}
||X^f||=\displaystyle\inf_{\substack{E\subseteq\Omega\\ \mathbb{P}(E)=0}}\left\{\sup_{\omega\in\Omega\backslash E}|X^f(\omega)|\right\}
\end{equation}
for $X^f=m_f-\hat{\mu}_f\in\mathscr{X}$.

Hence, the distance $d$ we use for $\mathscr{X}=\{m_f-\hat{\mu}_f,\,f\in\Psi\}$ is 
\begin{align}\label{distance}
d\,(X^f,X^{f'})=||X^f-X^{f'}||
\end{align}
for $X^f=m_f-\hat{\mu}_f\in\mathscr{X}$ and $X^{f'}=m_{f'}-\hat{\mu}_{f'}\in\mathscr{X}.$

Now, we present some definitions to introduce our theorems.

\begin{definition}[Entropy with bracketing]
$\mathscr{D}=\{D^{\lambda}: \lambda \in \boldsymbol{\Lambda}\}$ is a family of functions indexed by $\boldsymbol{\Lambda}$ and it is assumed that $||D^{\lambda}||<\infty$ for all $D^{\lambda}\in\mathscr{D}$. Given $\delta>0$, let $\{[{D_{j}^{\lambda}}^L,{D_{j}^{\lambda}}^U]\}_{j=1}^m\subseteq\mathscr{D}\times\mathscr{D}$ be a collection of pairs such that for each $\lambda\in\boldsymbol{\Lambda }$ there exists a $j=j(\lambda)\in\{1,2,\cdots,m\}$ such that:
\begin{enumerate}
\item ${D_{j}^{\lambda}}^L\leq D^{\lambda}\leq{D_{j}^{\lambda}}^U$ ;
\item $\resizebox{!}{!}{\textit{d}}\,({D_{j}^{\lambda}}^L,{D_{j}^{\lambda}}^U)\leq\delta.$
\end{enumerate}
Let $N(\delta)$ be the smallest value of \,$m$\, for which such a bracketing set $\{[{D_{j}^{\lambda}}^L,{D_{j}^{\lambda}}^U]\}_{j=1}^m$ exists. $N(\delta)$ is defined as the covering number of $\mathscr{D}$. $H(\delta)=\log N(\delta)$ is defined as the entropy of $\mathscr{D}$.
\end{definition}
Under the assumption of the so-called exponential class, Chen and Wu \cite{wwb} obtained the concentration inequalities for some empirical process of certain linear time series. Our result is obtained in a similar setting. The definition of exponential class is given first.
\begin{definition}[Exponential class]
    $\mathscr{D}=\{D^\lambda: \lambda \in \boldsymbol{\Lambda}\}$ is a family of functions indexed by $\boldsymbol{\Lambda}$ and it is assumed that $||D^{\lambda}||<\infty$ for all $D^{\lambda}\in\mathscr{D}$. For some constants $A, B, r>0$, the covering number $N(\kappa)\leq A\exp{\left(B\kappa^{-r}\right)}$ for all $0<\kappa\leq n$. $\mathscr{D}$ is called an exponential class.
\end{definition}

\begin{remark}
    The definition above quantifies the magnitudes of the class $\mathscr{D}$. It can be checked that it is a subclass of Donsker class (see van der Vaart and Wellner \cite{vw} for details). 
\end{remark}

Depending on the specific form of the exponential class, we have the following two theorems.

\begin{theorem}\label{thm1}
Let $\mathscr{X}=\{m_f-\hat{\mu}_f: f \in \Psi\}$ be a family of deviations indexed by $\Psi$. Assume that $\displaystyle\sup_{f\in\Psi} Var[f(X)]\leq \sigma^2$ for some $\sigma^2>0$. 
For $\varepsilon\in(0,1)$, if the following two conditions: 
\begin{itemize}
    \item $\mathscr{X}$ is an exponential class satisfying $N(\delta)\leq C\exp{\left(M\delta^{-2}\right)}$ for all $0<\delta\leq n$ and some constants $C, M>0$,
    \item for some constants $C_1$ and $C_2$ depending on $\sigma^2$,
\begin{equation}\label{entropy}
    \tfrac{n^{3/2}\varepsilon^2}{C_1}\geq \left\{\sqrt{M}\log{\left(\frac{n}{\frac{n^2\varepsilon}{C_2}\wedge\frac{n}{4}}\right)}+\sqrt{\log C}\left(n-\frac{n^2\varepsilon}{C_2}\wedge\frac{n}{4}\right)\mathds{1}_{\{C>1\}}\right\}\vee n^{1/2},
\end{equation}
\end{itemize}
hold, then we have
\begin{equation}\label{conclusion}
  \mathbb{P}\left(\displaystyle\sup_{f\in\Psi}\left|m_f-\hat{\mu}_f\right|\geq n^2\varepsilon\right)\leq C_3\exp\left(-\frac{n\varepsilon^2}{C_4}\right)
\end{equation}
for some constants $C_3$ and $C_4$ depending on $\sigma^2$.

\end{theorem}

The following theorem gives uniform exponential inequality for a wider family of exponential classes.

\begin{theorem}\label{thm2}
Let $\mathscr{X}=\{m_f-\hat{\mu}_f: f \in \Psi\}$ be a family of deviations indexed by $\Psi$. Assume that $\displaystyle\sup_{f\in\Psi} Var[f(X)]\leq \sigma^2$ for some $\sigma^2>0$. For $\varepsilon\in(0,1)$, if the following two conditions:

\begin{itemize}
    \item $\mathscr{X}$ is an exponential class satisfying $N(\delta)\leq C\exp{\left(M\delta^{-p}\right)}$ for all $0<\delta\leq n$ and some constants $C, M, p>0$, $p\neq2$,
    \item for some constants $C_1$ and $C_2$ depending on $\sigma^2$, \begin{equation}\label{entropy2}
    \tfrac{n^{3/2}\varepsilon^2}{C_1}\geq \left\{\frac{2\sqrt{M}}{2-p}\left(n^{1-\frac{p}{2}}-\left(\frac{n^2\varepsilon}{C_2}\wedge\frac{n}{4}\right)^{1-\frac{p}{2}}\right)+\sqrt{\log C}\left(n-\frac{n^2\varepsilon}{C_2}\wedge\frac{n}{4}\right)\mathds{1}_{\{C>1\}}\right\}\vee n^{1/2},
\end{equation}
\end{itemize}
hold, then we have
\begin{equation}\label{conclusion2}
  \mathbb{P}\left(\displaystyle\sup_{f\in\Psi}\left|m_f-\hat{\mu}_f\right|\geq n^2\varepsilon\right)\leq C_3\exp\left(-\frac{n\varepsilon^2}{C_4}\right)
\end{equation}
for some constants $C_3$ and $C_4$ depending on $\sigma^2$.
\end{theorem}

\begin{remark}
    When $\Psi$ is uncountable, $\mathbb{P}$ is understood as the outer probability (see, e.g., Van de Vaart \cite{vw}).
\end{remark}

\begin{remark}
    When $\Psi$ is a finite class of cardinality, say $|\Psi|=N$, Catoni's estimators satisfy (see Brownless et al. \cite{b} for more details)
$$\mathbb{P}\left(\displaystyle\sup_{f\in\Psi}|m_f-\hat{\mu}_f|\geq\sqrt{\frac{2\sigma^2\log(N/\delta)}{n(1-(2/n)\log(N/\delta))}}\right)\leq2\delta.$$
\end{remark}

\begin{remark}
The proof of Theorem \ref{thm1} and Theorem \ref{thm2} are identical except little difference in (\ref{difference}) below. We only present proof of Theorem \ref{thm1} for simplicity.
\end{remark}

\begin{remark}
For the proof, we employ a chaining technique with adaptive truncation similar to that of van de Geer \cite{va}. More generally, van de Geer et al. \cite{vaa} proposed a method called \textit{chaining along a tree} based on the Bernstein-Orlicz norm proposed in \cite{vaa}, and they demonstrated that entropy with bracketing allowed one to construct a finite tree chain to use chaining along a tree method. \end{remark}

The rest of the paper is organized as follows. The proof of the Theorem \ref{thm1} will be presented in Section 2. Section 3 describes an application to empirical risk minimization when the losses are not necessarily bounded and may have a distribution with heavy tails. Brownless et al. \cite{b} considered empirical risk minimization based on Catoni's estimator and developed performance bound using Talagrand's \textit{generic chaining} method (see, e.g., Talagrand \cite{t}). However, the upper bound they get depends on two distances, making that result complicated. On the contrary, our result is shown to be much more user-friendly, which could be regarded as an extension of the work of Brownless et al.\cite{b}.

\section{Proof of Theorem \ref{thm1}}

For convenience, in the proof, $c_i, i=1,\cdots, 16,$ will be constants depending on $C_1, C_2, C_3, C_4$ and $\sigma^2$. Denote $m_f-\hat{\mu}_f$ by $X^f$ for simplicity.

The following inequality will be useful in our proof.
By (\ref{deviations}), it is easily checked that
\begin{equation}\label{increment}
    \mathbb{P}\left(|X^f-X^{f'}|\geq x\right)\leq4\,\exp\left(-\frac{nx^2}{2(x^2+4\sigma^2)}\right)\,\quad(x>0),
\end{equation}
for $f,f'\in \Psi$, by triangle inequality.

Let 
$$I=\min\left\{i\geq1,2^{-i}\leq\frac{n\varepsilon}{2^3}\right\}$$
and write $H_i=H(2^{-i}n),\,i=0,\cdots, I$. 

By choice of $I$, we get that
\begin{equation}
    \frac{2^3}{n\varepsilon}\leq2^I<\frac{2^4}{n\varepsilon}.
\end{equation}

It is easily justified that, with proper choice of $C_2$, 
\begin{align}\label{difference}
\displaystyle\sum_{i=0}^I2^{-i}H_i^{1/2}&\leq\frac{2}{n}\int_{\tfrac{n}{2^{I+1}}}^n\sqrt{H(x)}dx\leq\frac{2}{n}\int_{\tfrac{n^2\varepsilon }{C_2}\wedge\frac{n}{4}}^n\sqrt{H(x)}dx\nonumber&\\
&= \frac{2}{n}\int_{\tfrac{n^2\varepsilon }{C_2}\wedge\frac{n}{4}}^n\sqrt{\log N(x)}dx\nonumber&\\
&\leq\frac{2}{n}\int_{\tfrac{n^2\varepsilon}{C_2}\wedge\frac{n}{4}}^n\sqrt{\log C+Mx^{-2}}dx\nonumber&\\
&\leq\frac{2}{n}\int_{\tfrac{n^2\varepsilon }{C_2}\wedge\frac{n}{4}}^n\sqrt{M}x^{-1}+\sqrt{\log C}\mathds{1}_{\{C>1\}}dx&\\
&=\frac{2}{n}\left\{\sqrt{M}\log{\left(\frac{n}{\frac{n^2\varepsilon}{C_2}\wedge\frac{n}{4}}\right)}+\sqrt{\log C}\left(n-\frac{n^2\varepsilon}{C_2}\wedge\frac{n}{4}\right)\mathds{1}_{\{C>1\}}\right\}\nonumber&\\
&\leq\frac{2n^{3/2}\varepsilon^2}{C_1}.\nonumber&
\end{align}

Since $H_i\geq0,\,i=0,\cdots, I$, we get
\begin{equation}
    \displaystyle\sum_{i=1}^{I}2^{-i}\left(\sum_{k=0}^{i}H_k\right)^{1/2}\leq\displaystyle\sum_{i=1}^{I}2^{-i}\left(\sum_{k=0}^{i}H_k^{1/2}\right)\leq\displaystyle\sum_{i=0}^{I}\left(\sum_{j=i}^{I}2^{-j}\right)H_i^{1/2}\leq2\displaystyle\sum_{i=0}^I2^{-i}H_i^{1/2}.
\end{equation}
Then, 
\begin{equation}
    \displaystyle\sum_{i=1}^I\left(\sum_{k=0}^i H_k\right)^{1/2}\leq2^I\displaystyle\sum_{i=1}^I 2^{-i}\left(\sum_{k=0}^{i}H_k\right)^{1/2}\leq2^{I+1}\displaystyle\sum_{i=0}^I2^{-i}H_i^{1/2}\leq\frac{2^6n^{1/2}\varepsilon}{C_1}.\label{ineq1}
\end{equation}
Set
\begin{equation}
\eta_i=\max\left\{\frac{C_1\left(\displaystyle\sum_{k=0}^i H_k\right)^{1/2}}{2^7 n^{1/2}\varepsilon}, \frac{1}{20}\sqrt{\frac{i}{2^i}}\right\},\quad i=1,\cdots, I.
\end{equation}

Furthermore, by (\ref{ineq1}), we obtain that
$$\displaystyle\sum_{i=1}^I\eta_i\leq1.$$

For $i=0,\cdots, I, \{[\tilde{X}_i^{f^L},\tilde{X}_i^{f^U}]\}$ is a $(2^{-i}n)$-bracketing set for $\mathscr{X}$. The subscript refers to the bracketing set rather than the member of a bracketing set. Denote the $N(2^{-i}n)$ by $N_i$
and then $\log N_i=H_i$. At the same time, for any given $X^f\in\mathscr{X}$, there is a pair $[\tilde{X}_i^{f^L},\tilde{X}_i^{f^U}]$ such that $$\tilde{X}_i^{f^L}\leq X^f\leq\tilde{X}_i^{f^U}\quad \text{and}\quad \quad \resizebox{!}{!}{\textit{d}}(\tilde{X}_i^{f^L},\tilde{X}_i^{f^U})\leq2^{-i}n.$$
Define
\begin{equation}
    X^{f^U}_i=\displaystyle\min_{k\leq i}\tilde{X}_k^{f^U},\quad X^{f^L}_i=\displaystyle\max_{k\leq i}\tilde{X}_k^{f^L},\quad \Delta_i^f=X^{f^U}_i-X^{f^L}_i,\quad i=0,\cdots, I,
\end{equation}
and
\begin{equation}
    v=\displaystyle\min\left\{i\geq0,\,\Delta_i^f\geq K_i\right\}\wedge I,\quad\quad K_i=\frac{2^{3-2i}}{\varepsilon\eta_{i+1}},\quad i=0,\cdots, I-1.
\end{equation}

It is worth mentioning that the pairs $[\tilde{X}_i^{f^L},\tilde{X}_i^{f^U}]$, $X^{f^U}_i$, $X^{f^L}_i$ and $\Delta_i^f,\,i=0,\cdots, I,$ as well as $v$ depend on $X^f\in\mathscr{X}$, although this is not expressed explicitly in our notations.

We have the following identity:
$$X^f=X^{f^L}_0+\displaystyle\sum_{i=0}^I\left(X^f-X^{f^L}_i\right)\mathds{1}_{\{v=i\}}+\displaystyle\sum_{i=1}^I\left(X^{f^L}_i-X^{f^L}_{i-1}\right)\mathds{1}_{\{v\geq i\}}.$$

Hence,
\begin{align} \label{all}
\mathbb{P}\left(\displaystyle\sup_{f\in\Psi}\left|X^f\right|\geq n^2\varepsilon\right)
&\leq\mathbb{P}\left(\displaystyle\sup_{f\in\Psi}\left|X^{f^L}_0\right|\geq\frac{n^2\varepsilon}{2}\right)&\nonumber\\
&+\mathbb{P}\left(\displaystyle\sup_{f\in\Psi}\left|\sum_{i=0}^I\left(X^f-X^{f^L}_i\right)\mathds{1}_{\{v=i\}}\right|\geq\frac{n^2\varepsilon}{4}\right)&\nonumber\\
&+\mathbb{P}\left(\displaystyle\sup_{f\in\Psi}\left|\sum_{i=1}^I\left(X^{f^L}_i-X^{f^L}_{i-1}\right)\mathds{1}_{\{v\geq i\}}\right|\geq\frac{n^2\varepsilon}{4}\right)&\\
&=\mathbb{P}_{\uppercase\expandafter{\romannumeral1}}+\mathbb{P}_{\uppercase\expandafter{\romannumeral2}}+\mathbb{P}_{\uppercase\expandafter{\romannumeral3}}.&\nonumber
\end{align}

As for $\mathbb{P}_{\uppercase\expandafter{\romannumeral1}}$,  by $(\ref{deviations})$, we have
\begin{align}\label{p11}    \mathbb{P}_{\uppercase\expandafter{\romannumeral1}}&\leq2N_0\exp\left(-\frac{n\left(\frac{n^2\varepsilon}{2}\right)^2}{2\left(\sigma^2+\left(\frac{n^2\varepsilon}{2}\right)^2\right)}\right)
\leq2\exp\left(H_0-\frac{n\left(\frac{n^2\varepsilon}{2}\right)^2}{2\left(\sigma^2+\left(\frac{n^2\varepsilon}{2}\right)^2\right)}\right)&\nonumber\\
&\leq2\exp\left(H_0-\frac{n^5\varepsilon^2}{2(n^4+4\sigma^2)}\right)\leq2\exp\left(H_0-\frac{n\varepsilon^2}{c_1}\right).
\end{align}

By condition (\ref{entropy}), we get 
\begin{align*}
    \frac{3n}{4}H_0^{1/2}&\leq\int_{\frac{n^2\varepsilon }{C_2}\wedge \frac{n}{4}}^n\sqrt{H(x)}dx\\
    &\leq\left\{\sqrt{M}\log{\left(\frac{n}{\frac{n^2\varepsilon}{C_2}\wedge\frac{n}{4}}\right)}+\sqrt{\log C}\left(n-\frac{n^2\varepsilon}{C_2}\wedge\frac{n}{4}\right)\mathds{1}_{\{C>1\}}\right\}\\
    &\leq\frac{n^{3/2}\varepsilon^2}{C_1},
\end{align*}
namely,
\begin{equation}\label{cond.}
    H_0\leq\frac{16n\varepsilon^4}{9C_1^2}\leq\frac{c_2n\varepsilon^2}{C_1^2}.
\end{equation}

Therefore, for $C_1$ sufficiently large,
\begin{equation}\label{p12}
    \mathbb{P}_{\uppercase\expandafter{\romannumeral1}}\leq2\exp\left(-\frac{n\varepsilon^2}{c_3}\right).
\end{equation}

For $i=0,\cdots, I-1$,
\begin{equation}\label{sum1}
    \Delta_i^f\mathds{1}_{\{v=i\}}\leq
    \frac{(\Delta_i^{f})^2}{K_i}\mathds{1}_{\{v=i\}}\leq\frac{\resizebox{!}{!}{\textit{d}}^2\left(X^{f^L}_i, X^{f^U}_i\right)}{K_i}\leq\frac{(2^{-i}n)^2}{K_i}\leq\frac{n^2\varepsilon\eta_{i+1}}{2^3}.
\end{equation}

For $i=I$,
\begin{equation}\label{sum2}
    \Delta^f_I\mathds{1}_{\{v=I\}}\leq d\left(X^{f^L}_I,X^{f^U}_I\right)\leq\frac{n}{2^I}\leq\frac{n^2\varepsilon}{2^3}.
\end{equation}

Since $\displaystyle\sum_{i=1}^I\eta_i\leq1$, it follows from (\ref{sum1}) and (\ref{sum2}) that
$$\displaystyle\sum_{i=0}^I\Delta^f_i\mathds{1}_{\{v=i\}}\leq\sum_{i=0}^{I-1}\frac{n^2\varepsilon\eta_{i+1}}{2^3}+\frac{n^2\varepsilon}{2^3}\leq\frac{n^2\varepsilon}{4}.$$

This implies that 
\begin{align*}
    \mathbb{P}_{\uppercase\expandafter{\romannumeral2}}&\leq\mathbb{P}\left(\displaystyle\sup_{f\in\Psi}\sum_{i=0}^I\left|\Delta_i^f\mathds{1}_{\{v=i\}}\right|\geq\frac{n^2\varepsilon}{4}\right)&\\
&\leq\mathbb{P}\left(\displaystyle\sup_{f\in\Psi}\left|\Delta_0^f\mathds{1}_{\{v=0\}}\right|\geq\frac{n^2\varepsilon}{2^3}\right)&\\
&+\mathbb{P}\left(\displaystyle\sup_{f\in\Psi}\sum_{i=1}^I\left|\Delta_i^f\mathds{1}_{\{v=i\}}\right|\geq\frac{n^2\varepsilon}{2^3}\right)&\\
&=\mathbb{P}_{\uppercase\expandafter{\romannumeral2}}^{(1)}+\mathbb{P}_{\uppercase\expandafter{\romannumeral2}}^{(2)}.&
\end{align*}

For $\mathbb{P}_{\uppercase\expandafter{\romannumeral2}}^{(1)}$, similar to (\ref{p11}), (\ref{cond.}) and (\ref{p12}), we obtain
\begin{align}\label{p21}
\mathbb{P}_{\uppercase\expandafter{\romannumeral2}}^{(1)}&\leq2N_0\exp\left(-\frac{n\left(\frac{n^2\varepsilon}{2^3}\right)^2}{2\left(\sigma^2+\left(\frac{n^2\varepsilon}{2^3}\right)^2\right)}\right)\nonumber&\\
&\leq2\exp\left(H_0-\frac{n\left(\frac{n^2\varepsilon}{2^3}\right)^2}{2\left(\sigma^2+\left(\frac{n^2\varepsilon}{2^3}\right)^2\right)}\right)&\\
&\leq2\exp\left(H_0-\frac{n\varepsilon^2}{c_4}\right)\leq2\exp\left(-\frac{n\varepsilon^2}{c_5}\right).&\nonumber
\end{align}

The following observation will be useful when dealing with $\mathbb{P}_{\uppercase\expandafter{\romannumeral2}}^{(2)}$,
\begin{equation*}\label{observ}
\left(\displaystyle\sum_{k=0}^iH_k\right)^{1/2}\leq\frac{c_6n^{1/2}\varepsilon\eta_i}{C_1}.
\end{equation*}

As a consequence, combining with the definition of $\eta_i$, condition (\ref{increment}) and condition (\ref{entropy}), for $C_1$ sufficiently large, we have
\begin{align}\label{p22}
    \mathbb{P}_{\uppercase\expandafter{\romannumeral2}}^{(2)}&\leq\displaystyle\sum_{i=1}^I\mathbb{P}\left(\displaystyle\sup_{f\in\Psi}\left|\Delta_i^f\mathds{1}_{\{v=i\}}\right|\geq\frac{n^2\varepsilon\eta_i}{2^3}\right)\nonumber&\\
      &\leq\sum_{i=1}^I4N_i\exp\left(-\frac{n\left(\frac{n^2\varepsilon\eta_i}{2^3}\right)^2}{2\left(\left(\frac{n^2\varepsilon\eta_i}{2^3}\right)^2+4\sigma^2\right)}\right)\nonumber&\\
      &\leq\sum_{i=1}^I4\exp\left(H_i-\frac{n\left(\frac{n^2\varepsilon\eta_i}{2^3}\right)^2}{2\left(\left(\frac{n^2\varepsilon\eta_i}{2^3}\right)^2+4\sigma^2\right)}\right)\nonumber&\\
      &\leq\displaystyle\sum_{i=1}^I4\exp\left(\sum_{k=0}^iH_k-\frac{n\left(\frac{n^2\varepsilon\eta_i}{2^3}\right)^2}{2\left(\left(\frac{n^2\varepsilon\eta_i}{2^3}\right)^2+4\sigma^2\right)}\right)&\\
    &\leq\displaystyle\sum_{i=1}^I4\exp\left(\sum_{k=0}^iH_k-\frac{n\varepsilon^2\eta_i^2}{c_7}\right)\nonumber&\\
    &\leq\displaystyle\sum_{i=1}^I4\exp\left(-\frac{n\varepsilon^2\eta_i^2}{c_8}\right)\nonumber\leq\displaystyle\sum_{i=1}^I4\exp\left(-\frac{n\varepsilon^2 i}{c_92^i}\right)\nonumber&\\
    &\leq\displaystyle\sum_{i=1}^I4\exp\left(-\frac{n\varepsilon^2 i}{c_9 2^I}\right)\nonumber\leq\displaystyle\sum_{i=1}^I4\exp\left(-\frac{C_1n\varepsilon i}{c_{10}}\right)\leq c_{11}\exp\left(-\frac{n\varepsilon^2}{c_{12}}\right).&\\
\end{align}

As for $\mathbb{P}_{\uppercase\expandafter{\romannumeral3}}$, we may write
\begin{align*}
    \mathbb{P}_{\uppercase\expandafter{\romannumeral3}}&\leq\mathbb{P}\left(\displaystyle\sup_{f\in\Psi}\left|\left(X^{f^L}_1-X^{f^L}_0\right)\mathds{1}_{\{v\geq 1\}}\right|\geq\frac{n^2\varepsilon}{2^3}\right)&\\
&+\mathbb{P}\left(\displaystyle\sup_{f\in\Psi}\left|\sum_{i=2}^I\left(X^{f^L}_i-X^{f^L}_{i-1}\right)\mathds{1}_{\{v\geq i\}}\right|\geq\frac{n^2\varepsilon}{2^3}\right)&\\
&=\mathbb{P}_{\uppercase\expandafter{\romannumeral3}}^{(1)}+\mathbb{P}_{\uppercase\expandafter{\romannumeral3}}^{(2)}.&
\end{align*}

It is noteworthy that 
$\left|X^{f^L}_i-X^{f^L}_{i-1}\right|\leq\Delta^f_{i-1},\, i=1,\cdots, I$, which allows us to deal with $\mathbb{P}_{\uppercase\expandafter{\romannumeral3}}^{(1)}$ and $\mathbb{P}_{\uppercase\expandafter{\romannumeral3}}^{(2)}$ in a similar way as $\mathbb{P}_{\uppercase\expandafter{\romannumeral2}}^{(1)}$ and $\mathbb{P}_{\uppercase\expandafter{\romannumeral2}}^{(2)}$.

Similar to (\ref{p21}), we obtain 
\begin{align}\label{p31}
\mathbb{P}_{\uppercase\expandafter{\romannumeral3}}^{(1)}&\leq\mathbb{P}\left(\displaystyle\sup_{f\in\Psi}\left|\Delta_0^f\mathds{1}_{\{v=0\}}\right|\geq\frac{n^2\varepsilon}{2^3}\right)\nonumber&\\
&\leq2N_0\exp\left(-\frac{n\left(\frac{n^2\varepsilon}{2^3}\right)^2}{2\left(\sigma^2+\left(\frac{n^2\varepsilon}{2^3}\right)^2\right)}\right)\nonumber&\\
&\leq2\exp\left(H_0-\frac{n\left(\frac{n^2\varepsilon}{2^3}\right)^2}{2\left(\sigma^2+\left(\frac{n^2\varepsilon}{2^3}\right)^2\right)}\right)&\\
&\leq2\exp\left(H_0-\frac{n\varepsilon^2}{c_{13}}\right)\leq2\exp\left(-\frac{n\varepsilon^2}{c_{14}}\right).&\nonumber
\end{align}

By similar argument as in (\ref{p22}), we get
\begin{align}\label{p32}
    \mathbb{P}_{\uppercase\expandafter{\romannumeral3}}^{(2)}&\leq\mathbb{P}\left(\displaystyle\sup_{f\in\Psi}\sum_{i=1}^{I-1}\left|\Delta_i^f\mathds{1}_{\{v=i\}}\right|\geq\frac{n^2\varepsilon}{2^3}\right)\nonumber&\\
    &\leq\displaystyle\sum_{i=1}^{I-1}\mathbb{P}\left(\displaystyle\sup_{f\in\Psi}\left|\Delta_i^f\mathds{1}_{\{v=i\}}\right|\geq\frac{n^2\varepsilon\eta_i}{2^3}\right)&\\
    &\leq c_{15}\exp\left(-\frac{n\varepsilon^2}{c_{16}}\right).&\nonumber
\end{align}

Combining (\ref{all}), (\ref{p12}), (\ref{p21}), (\ref{p22}), (\ref{p31}) and (\ref{p32}), we complete the proof of the Theorem $\ref{thm1}$.

\section{Application}
In this section, we describe the application of Theorem \ref{thm1} and Theorem \ref{thm2} to empirical risk minimization for heavy-tailed losses.

For a random feature variable $Z$ taking values in some measurable space $\mathcal{Z}$ and a random target variable $Y$ taking values in $\mathbb{R}$, let $\Psi$ be a set of real-valued functions defined on $\mathcal{Z}$. It is anticipated to find a \textit{prediction function} $g\in\Psi$ that describes the relationship between $Z$ and $Y$, which is assumed to follow a joint probability $P(Z, Y)$. To this end, a \textit{loss function} $\ell:\mathbb{R}\times\mathbb{R}\to\mathbb{R}_+$ is defined to penalize the differences between prediction functions $g(z)$ and actual targets $y$. Then, minimize the expectation of the loss function $\ell$ over the data distribution $\mathbb{P}$, also known as the \textit{expected risk}:
$$R(g)=\mathbb{E}\ell\left(g(Z), Y\right)=\int\ell\left(g(z),y\right) d\mathbb{P}(z,y).$$
Throughout this section, it is assumed that $\displaystyle\inf_{g\in\Psi}R(g)$ is attained for some $g^*\in\Psi$ and denote the optimal risk $\displaystyle\inf_{g\in\Psi}R(g)$ by $m^*$. For example, in regression, $(Z, Y)\in\mathbb{R}^d\times\mathbb{R}$, $\ell(g(Z),Y)=(Y-g(Z))^2$ for $g$ in some class $\mathscr{G}$, such as the class of linear functions. In this setting, $g^*(z)=\mathbb{E}[Y|Z=z]$ is the conditional expectation.

Unfortunately, the distribution $\mathbb{P}$ is unknown in most practical situations. Instead, a set of ``training data" $(z_1, y_1),\cdots, (z_n,y_n)$, with $(z_i,y_i)\sim\mathbb{P}$ for all $i=1,\cdots, n$ are usually available. Using the training data, the distribution $\mathbb{P}$ might be approximated by the \textit{empirical distribution} 
$$\mathbb{P}_{\delta}(z,y)=\frac{1}{n}\displaystyle\sum_{i=1}^n\delta(z=z_i, y=y_i),$$
where $\delta(z=z_i, y=y_i)$ is a Dirac mass centered at $(z_i, y_i)$. Using the empirical distribution $\mathbb{P}_{\delta}$, the expected risk can now be approximated by the \textit{empirical risk}:
\begin{equation}\label{emp risk}
    R_{\delta}(g)=\int\ell\left(g(z),y\right)d\mathbb{P}_{\delta}(z,y)=\frac{1}{n}\displaystyle\sum_{i=1}^n\ell\left(g(z_i), y_i\right).
\end{equation}

The \textit{empirical risk minimizer} $g_{ERM}$ is defined as
$$g_{ERM}=\mathop{\arg\min}\limits_{g\in\Psi}R_{\delta}(g). $$

The quality of empirical risk minimization is measured by the \textit{excessive risk} $$R_{\delta}(g_{ERM})-m^*.$$ 
In this section, we mainly focus on excessive risk.

Learning the function $g$ by minimizing (\ref{emp risk}) is known as the \textit{empirical risk minimization} principle (see, e.g., Vapnik \cite{v}). 

In this section, to simplify the notation, $X_i$ will denote the pair $(Z_i,Y_i)$, the function $f$ substitutes $\ell\left(g(\cdot),\cdot\right)$ and $m_f$ represents $R(g)$. It might be assumed that $f\in\mathcal{F}$ and $X$ is independent from $(X_i)_{i=1}^{n}$ with the same marginal distribution $F$. Also, it might be supposed that $X$ takes values in some measurable space $\mathcal{S}$ and $\mathcal{F}$ is a set of functions defined on $\mathcal{S}$.

When $f\in\mathcal{F}$ are all uniformly bounded or the random variables $f(X)$ have sub-Gaussian tails for all $f\in\mathcal{F}$, the performance of empirical risk minimization has been well understood, see Adamczak \cite{ada}, Koltchinskii \cite{kol} and van der Vaart and Wellner \cite{vw}. However, when the functions $f$ are no longer uniformly bounded and the random variables $f(X)$ may have a heavy tail, empirical risk minimization may behave poorly since the empirical mean has become a bad estimator for the expected value. Lugosi and Mendelson \cite{l3} showed why empirical risk minimization failed in the mean-squared sense. This motivates us to consider a robust version of empirical risk minimization based on minimizing Catoni's estimators.

For fixed $f\in\mathcal{F}$, Catoni's estimator is denoted by $\hat{\mu}_f$. It is intuitive to define an empirical risk minimizer from the class $\mathcal{F}$ that minimizes Catoni's estimators, which is also how Brownless et al. \cite{b} defined the empirical risk minimizer. So define 
\begin{equation}\label{cat.erm}
   \hat{f}=\mathop{\arg\min}\limits_{f\in\mathcal{F}}\hat{\mu}_f,
\end{equation}
and the performance of the empirical risk minimization is evaluated by $$m_{\hat{f}}-m^*.$$

Some concepts will be introduced to describe the results in Brownless et al. \cite{b}.

For the class $\mathcal{F}$, the distance $d_1$ is defined, for $f, f'\in\mathcal{F}$ by
$$d_1(f, f')=\left(\mathbb{E}\left[\left(f(X)-f'(X)\right)^2\right]\right)^{1/2},$$
and the distance $d_2$ is 
$$d_2(f, f')=\displaystyle\sup_{x\in\mathcal{S}}\left|f(x)-f'(x)\right|.$$

Let $T$ be a metric space. An increasing sequence $(\mathcal{A}_n)$ of partitions of $T$ is called \textit{admissible} if $\#\mathcal{A}_n\leq N_n=2^{2^n}$ for all $n=0, 1, 2,\cdots$. For any $t\in T$, denote by $A_n(t)$ the unique element of $\mathcal{A}_n$ that contains $t$. Let $\Delta(A)$ denote the diameter of the set $A\subseteq T$. Define, for $\beta=1,2,$
\begin{equation}\label{acc}
    \gamma_{\beta}(T,d)=\displaystyle\inf_{\mathcal{A}_n}\sup_{t\in T}\sum_{n\geq0}2^{n/\beta}\Delta(A_n(t)),
\end{equation}
where the infimum is taken over all admissible sequences. 

What Brownless et al. \cite{b} obtained might be rewritten as the following theorem.
\begin{theorem}\label{b1}
Consider the setup as above and assume that $\displaystyle\sup_{f\in\mathcal{F}}Var\left(f(X)\right)\leq \sigma^2$ for some $\sigma^2>0$. Let $\varepsilon\in(\log6,+\infty)$. Suppose that $\hat{f}$ is selected from $\mathcal{F}$ by minimizing Catoni's estimator with parameter $\alpha$. Then there exists a constant $L$ such that, under the condition 
$$6\left(\alpha v+\frac{2\varepsilon}{n\alpha}\right)+L\varepsilon\left(\frac{\gamma_2(\mathcal{F},d_1)}{\sqrt{n}}+\frac{\gamma_1(\mathcal{F},d_2)}{n}\right)\leq\frac{1}{\alpha},$$
we have
\begin{equation}\label{twodist}
\mathbb{P}\left(m_{\hat{f}}-m^*\geq 6\left(\alpha v+\frac{2\varepsilon}{n\alpha}\right)+L\varepsilon\left(\frac{\gamma_2(\mathcal{F},d_1)}{\sqrt{n}}+\frac{\gamma_1(\mathcal{F},d_2)}{n}\right)\right)\leq 6\exp(-\varepsilon).
\end{equation}
\end{theorem}

\begin{remark}
It can be seen from (\ref{twodist}) that two distances $d_1$ and $d_2$ are needed to obtain the upper bound for excessive risk in Theorem \ref{b1}.
\end{remark}

Only with one distance $d$ defined in (\ref{distance}), applying Theorem \ref{thm1} and Theorem \ref{thm2} gives rise to the following two theorems, repspectively.

\begin{theorem}\label{cor}
Consider the setup described above and assume that $\displaystyle\sup_{f\in\mathcal{F}}Var\left(f(X)\right)\leq \sigma^2$ for some $\sigma^2>0$. For $\varepsilon\in(0,1)$, if the following two conditions:
\begin{itemize}
\item $\mathscr{X}=\{m_f-\hat{\mu}_f,\,f\in\mathcal{F}\}$
is an exponential class with $N(\delta)\leq C\exp{\left(M\delta^{-2}\right)}$ for all $0<\delta\leq n$ and some constants $C, M>0$,
\item for some constants $C_1$ and $C_2$ depending on $\sigma^2$, $$\tfrac{n^{3/2}\varepsilon^2}{C_1}\geq \left\{\sqrt{M}\log{\left(\frac{n}{\frac{n^2\varepsilon}{C_2}\wedge\frac{n}{4}}\right)}+\sqrt{\log C}\left(n-\frac{n^2\varepsilon}{C_2}\wedge\frac{n}{4}\right)\mathds{1}_{\{C>1\}}\right\}\vee n^{1/2},$$
\end{itemize}
hold, then 
$$\mathbb{P}\left(m_{\hat{f}}-m^*\geq n^2\varepsilon\right)\leq C_3\exp\left(-\frac{n\varepsilon^2}{C_4}\right),$$
for some constants $C_3$ and $C_4$ depending on $\sigma^2$.

\end{theorem}

\begin{theorem}\label{cor2}
Consider the setup described above and assume that $\displaystyle\sup_{f\in\mathcal{F}}Var\left(f(X)\right)\leq v$ for some $v>0$. For $\varepsilon\in(0,1)$, if the following two conditions:
\begin{itemize}
    \item $\mathscr{X}=\{m_f-\hat{\mu}_f,\,f\in\mathcal{F}\}$
is an exponential class with $N(\delta)\leq C\exp{\left(M\delta^{-p}\right)}$ for all $0<\delta\leq n$ and some constants $C, M, p>0$, $p\neq2$,
 \item for some constants $C_1$ and $C_2$ depending on $\sigma^2$, $$\tfrac{n^{3/2}\varepsilon^2}{C_1}\geq \left\{\frac{2\sqrt{M}}{2-p}\left(n^{1-\frac{p}{2}}-\left(\frac{n^2\varepsilon}{C_2}\wedge\frac{n}{4}\right)^{1-\frac{p}{2}}\right)+\sqrt{\log C}\left(n-\frac{n^2\varepsilon}{C_2}\wedge\frac{n}{4}\right)\mathds{1}_{\{C>1\}}\right\}\vee n^{1/2},$$
\end{itemize}
hold, then 
$$\mathbb{P}\left(m_{\hat{f}}-m^*\geq n^2\varepsilon\right)\leq C_3\exp\left(-\frac{n\varepsilon^2}{C_4}\right),$$

\end{theorem}

\begin{remark}
The covering number in these two theorems is defined with respect to the distance $d$ in (\ref{distance}). We gave a bound for excessive risk by one distance. Comparatively, our theorems demonstrate a much more user-friendly performance than Theorem \ref{b1}.
\end{remark}

\begin{remark}
The robust version of empirical risk minimization is not unique. In (\ref{cat.erm}), we replace the empirical mean with Catoni's estimator in an empirical risk minimization framework. Other robust mean estimators might also be employed to get another robust version of empirical risk minimization. For example, Mathieu and Minsker \cite{m2} used Minsker's estimator \cite{m} for robust empirical risk minimization and gave high-confidence bounds for excessive risk (see Mathieu and Minsker \cite{m2} for details).
\end{remark}

{\bf Proof of Theorem \ref{cor} and Theorem \ref{cor2}.}   
    
Note that 
$$m_{\hat{f}}-m^*=(m_{\hat{f}}-\hat{\mu}_{\hat{f}})+(\hat{\mu}_{\hat{f}}-m^*)\leq2\displaystyle\sup_{f\in\mathcal{F}}|m_f-\hat{\mu}_f|,$$
and then we have
\begin{equation}\label{observe}
\mathbb{P}\left(m_{\hat{f}}-m^*\geq n^2\varepsilon\right)\leq\mathbb{P}\left(\displaystyle\sup_{f\in\mathcal{F}}|m_f-\hat{\mu}_{f}|\geq\frac{n^2\varepsilon}{2}\right).
\end{equation}
Under the conditions listed in the theorem, we know that 
\begin{equation}\label{observe2}
\mathbb{P}\left(\displaystyle\sup_{f\in\mathcal{F}}|m_f-\hat{\mu}_{f}|\geq\frac{n^2\varepsilon}{2}\right)\leq C_3\exp\left(-\frac{n\varepsilon^2}{C_4}\right),
\end{equation}
for some constants $C_3$ and $C_4$.

A combination of (\ref{observe}), (\ref{observe2}) and Theorem \ref{thm1} (or Theorem \ref{thm2}, respectively) results in the Theorem \ref{cor} (or Theorem \ref{cor2}, respectively). 

\vskip3mm

\section*{Acknowledgments}
This work  was supported by  the National Natural Science Foundation of China (Nos. 12071257 and 11971267 );   National Key R$\&$D Program of China (No. 2018YFA0703900); Shandong Provincial Natural Science Foundation (No. ZR2019ZD41); and Young Scholars Program of Shandong University.

\vskip3mm

\bibliographystyle{amsplain}


\end{document}